\documentclass{alf-j-l}
\usepackage{verbatim} 
\usepackage{syntonly}
\usepackage{url}  

\ifx\preambleloaded\relax
   \else\let\preambleloaded\relax\fi
\def\Z{{\mathbb Z}} \def\Q{{\mathbb Q}} 

\def\F{{\mathbb F}} 

\def\={\equiv}

\def\\{\cr} 

\def\acclap#1{\raise\hgtsig\hbox to0pt{$#1$\hss}}
\newdimen\hgtsig
\setbox0=\hbox{$\displaystyle{\sum}$}
\hgtsig=\ht0\relax
\advance\hgtsig by -1.75ex

\newbox\boxW\newdimen\dimW
\def\heighten#1{%
\setbox\boxW\hbox{$\displaystyle #1$}
\dimW=1.04\ht\boxW\advance\dimW by 1.00pt
\vbox to \dimW{}}

\DeclareMathAlphabet{\Bi}{OT1}{cmm}{b}{it}  

\providecommand{\bysame}{\makebox[3em]{\hrulefill}\thinspace}
\def\con#1=#2(#3){#1\equiv#2\pod{#3}}
                            {\end{enumerate}}
                            {\end{enumerate}}

\def\cfraci#1#2{#1_0+{
\let\ds\displaystyle
\def\vl{\ds1\vrule width0pt depth.5ex height2ex\over}
\vl{\ds#1_1 + {\vl{\ds#1_2 +{\vl\hskip.5em\ddots}}}}}}

\def\divides{{\mathchoice{\mathrel{\bigm|}}{\mathrel{\bigm|}}{\mathrel{|}}%
{\mathrel{|}}}}
\def\Div{\divides}
\def\notdivides{\mathrel{\kern-3pt\not\!\kern3.5pt\bigm|}}

\newbox\boxA
\newbox\boxB
\newdimen\dimA
\newdimen\dimB
\newdimen\dimC

\def\house#1{
\setbox\boxA\hbox{$\displaystyle #1$}
\dimA=1.04\wd\boxA\advance\dimA by 2pt
\dimB=1.04\ht\boxA\advance\dimB by 2pt
\dimC=0.05\wd\boxA
\hskip\dimC\hskip1pt\hbox to \dimA
{\vrule\vbox to \dimB{\hsize=\dimA
\hrule\vfill \centerline{\box\boxA}}\vrule}
\hskip 0.8pt\hskip1pt\hskip\dimC}

\def\cf {continued fraction}
\def\pq {partial quotient}

\def\cfe{continued fraction expansion}

\newcount\hours
\newcount\minutes
\def \SetTime{\hours=\time
\global\divide\hours by 60
\minutes=\hours
\multiply\minutes by 60
\advance\minutes by-\time
\global\multiply\minutes by-1 }
\SetTime
\def \now{\number\hours:\ifnum\minutes<10 0\fi\number\minutes}
\def \Now{\today\ $[$\now$]$}

\newif\ifMacTextures
\MacTexturestrue

\ifMacTextures

\gdef\EPSF#1by#2(#3){%
\vbox to #2{\hrule width #1 height 0pt depth 0pt%
\vfill\special{illustration #3}}}%

\gdef\scaledEPSF#1by#2(#3 scaled #4){{%
\dimen0=#1 \dimen1=#2%
\divide\dimen0 by 1000 \multiply\dimen0 by #4%
\divide\dimen1 by 1000 \multiply\dimen1 by #4%
\EPSF \dimen0 by \dimen1 (#3 scaled #4)}}%
\else
\input epsf

\gdef\EPSF#1by#2(#3){%
\vbox to #2{\hrule width #1 height 0pt depth 0pt\vfill \epsfbox{#3}}}%

\gdef\EPSF#1by#2(#3){\epsfbox{#3}}%

\gdef\scaledEPSF#1by#2(#3 scaled #4){{%
\dimen0=#1 \dimen1=#2%
\divide\dimen0 by 1000 \multiply\dimen0 by #4%
\divide\dimen1 by 1000 \multiply\dimen1 by #4%
\epsfxsize=\dimen0\epsfbox{#3}}}%
\fi
\MacTexturesfalse

\def\ceNTRelogo{\vtop{\baselineskip10truept
\hsize= 0.8125in \smallcentrelogo \logofont
\smallskip
\centerline{ceNTRe}
\centerline{Sydney, Australia 2071}} }

\def\ds{\displaystyle}

\theoremstyle{plain}  
\newtheorem{theorem}{Theorem}

\newtheorem{proposition}[theorem]{Proposition} 

\newtheorem*{theorem*}{Theorem}
\newtheorem{corollary}[theorem]{Corollary}

\theoremstyle{definition}

\newtheorem*{Remark*}{Remark}
\newtheorem*{Remarks*}{Remarks}

\newtheorem*{example*}{Example} 
\newtheorem*{guess*}{Guess}

\theoremstyle{remark}

\newtheorem*{remark*}{Remark} 
\newtheorem*{remarks*}{Remarks}

\newtheoremstyle{aside}
   {6pt}
   {6pt}
   {\footnotesize}
   {}
   {\scshape}
   {:}
   {.5em}
   {}

\theoremstyle{aside}

\newtheorem*{aside*}{Concluding Aside}

\def\e{{\smash{\overline{e}}}}

\PII{Centre for Number Theory Research, 1 Bimbil Place, Killara,
Australia 2071}%

\copyrightinfo{\number\year}{Alfred J van der Poorten AM and
Christine S. Swart}

\begin{document}
\renewcommand\thesubsection{\arabic{subsection}}

\def\currentvolume{175}
\def\currentissue{Draft of\ }
\def\paperdate{\today}
\def\ISSN{}
\pagespan{69}{\pageref{page:lastpage}}  

\def\ceNTRelogo{\relax}

\title {Recurrence Relations\rlap{\smash{\qquad\qquad\ceNTRelogo}}\\
for Elliptic Sequences\,:\\every Somos~$4$ is a Somos~$k$} 

\author{Alfred J. van der Poorten}
\address{Centre for Number Theory Research, 1 Bimbil Place, Killara, Sydney,
NSW 2071, Australia}
\email{alf@math.mq.edu.au (Alf van der Poorten)}
\author{Christine S. Swart}
\address{Department of Mathematics, Statistics, and Computer Science\\
The University of Illinois at Chicago\\
Chicago, IL 60607-7045 USA}
\email{swart@math.uic.edu (Christine Swart)}

\thanks{The first author was supported in part by a grant from the
Australian Research Council.}


\date{\Now.}


\begin{abstract} In his `Memoir on Elliptic Divisibility Sequences',
Morgan Ward's definition of the said sequences has the
remarkable feature that it does not become at all clear until deep
into the paper that there exist nontrivial such sequences. Even
then, Ward's proof of coherence of his definition relies on
displaying a sequence of values of quotients of Weierstra\ss\
$\sigma$-functions. We give a direct proof of coherence and show,
rather more generally, that a sequence defined by a so-called
Somos~relation of gap~$4$ always  also is given by a three term
Somos~relation of all larger gaps $5$,
$6$, $7$, $\ldots\,$. 
\end{abstract}

\maketitle
\pagestyle{myheadings}\markboth{{\headlinefont Alf van der
Poorten and Christine Swart}}{{\headlinefont Recurrence
relations for elliptic sequences}}


\subsection{Morgan Ward's elliptic sequences} In his `Memoir on
elliptic divisibility sequences'~\cite{Wa}, Morgan Ward in effect
(thus, for all practical purposes) defines anti\-symmetric
double-sided sequences
$(W_h)$, that is with $W_{-h}=-W_h$, by requiring that, for all
integers
$h$,
$m$, and $n$,
\begin{equation} \label{eq:redundant'}
W_{h-m}W_{h+m}W_n^2+W_{n-h}W_{n+h}W_m^2
+W_{m-n}W_{m+n}W_h^2=0\,.\tag{\ref{eq:redundant}$'$}
\end{equation} 
If one dislikes double-sided sequences then one rewrites
\eqref{eq:redundant'} less elegantly as 
\begin{equation} \label{eq:redundant}
W_{h-m}W_{h+m}W_n^2=W_{h-n}W_{h+n}W_m^2 - W_{m-n}W_{m+n}W_h^2\,,
\end{equation}
just for $h\ge m\ge n$. 
In any case, \eqref{eq:redundant} seems more dramatic than it is. An
easy exercise confirms that if $W_1=1$ then~\eqref{eq:redundant} is
equivalent to just 
\begin{equation} \label{eq:general}
W_{h-m}W_{h+m}=W_m^2 W_{h-1}W_{h+1} -W_{m-1}W_{m+1}W_h^2
\end{equation}
for all integers $h\ge m$. Indeed, \eqref{eq:general} is just a special
case of \eqref{eq:redundant}. However, given~\eqref{eq:general}, obvious
substitutions in
\eqref{eq:redundant} quickly show one may return from
\eqref{eq:general} to the apparently more general \eqref{eq:redundant}.

But there is a drama here.  The recurrence relation
$$
W_{h-2}W_{h+2}=W_2^2 W_{h-1}W_{h+1} -W_{1}W_{3}W_h^2\,,
$$ 
and 
non-zero initial values $W_1=1$, $W_2$, $W_3$, $W_4$, already
suffices to produce the complete sequence!  Thus \eqref{eq:general}
for all
$m$ is entailed by its special case m = 2.

We could show directly that the case $m=3$ follows, see a remark in
\cite{Sw}, or a footnote in the corresponding discussion in \cite{169}, but
the case $m=4$, if done asymetrically as in subsequent remarks
of~\cite{169}, plainly was not worth the effort.  Plan B, to look it
up, fared little better. In her thesis
\cite{Shi}, Rachel Shipsey shyly refers the reader back to Morgan
Ward's memoir \cite{Wa}; but at first glance Ward seems not to
comment on the matter at all, having \emph{defined} his sequences by
\eqref{eq:general}.  Of course, Ward does comment. The issue is
whether
\eqref{eq:general} is \emph{coherent}: do different $m$ yield the one
sequence? Ward notes that if $\sigma$ is the Weierstra\ss\
$\sigma$-function then a sequence $\bigl(\sigma(hu)/\sigma(u)^{h^2}\,\bigr)$
satisfies
\eqref{eq:general} for all $m$.  He then painfully shows there is a
related cubic curve for every choice of initial
values.  Whatever, a much more direct
argument would be much more satisfying: we supply such an argument
below.

\subsection{Somos sequences} Some years ago, Michael Somos,
see~\cite{Somos}, 
\emph{inter~alia} asked for the inner meaning of the behaviour of
the sequences $(C_h)=(\dots\,$,  $2$, $1$, $1$, $1$, $1$, $2$, $3$,
$7$, $23$, $59$,
$\ldots\,)$ defined by
$C_{h-2}C_{h+2}=C_{h-1}C_{h+1}+C_{h}^2$; and of $(B_h)= (\dots\,$,
$2$, $1$, $1$, $1$, $1$, $1$, $2$, $3$, $5$, $11$, $37$, $83$, 
$\ldots\,)$ defined by
$B_{h-2}B_{h+3}=B_{h-1}B_{h+2}+B_{h}B_{h+1}$: that is, the sequences
$4$-Somos and $5$-Somos \cite[A006720 and A006721]{Sl}. More
generally, of course, one may both vary the `initial' values and
coefficients and generalise the `gap' to
$2m$ or
$2m+1$ by studying Somos~$2m$, respectively Somos~$2m+1$, namely
sequences satisfying the respective recursions
$$
D_{h-m}D_{h+m}=\sum_{i=1}^m \kappa_iD_{h-m+i}D_{h+m-i} \text{\ or\ }
D_{h-m}D_{h+m+1}=\sum_{i=1}^m \kappa_iD_{h-m+i}D_{h+m-i}\,.
$$
Direct, but somewhat painful, attacks allow one to prove that in
fact a Somos~$4$ always is a Somos~$5$, Somos~$6$, and Somos~$8$. For
example, see \cite{169}, $4$-Somos satisfies all of
\begin{align*}
C_{h-2}C_{h+3}&=-C_{h-1}C_{h+2}+5C_{h}C_{h+1}\,,\\
C_{h-3}C_{h+3}&=C_{h-1}C_{h+1}+5C_{h}^2\,,\\
C_{h-4}C_{h+4}&=25C_{h-1}C_{h+1}-4C_{h}^2\,.
\end{align*}
In the light of such results one feels some confidence that in
general a Somos~$4$ is a Somos~$k$, for all $k=5$, $6$, $7$,
$\ldots\,$: indeed, it is that which we show below.
\subsection{Elliptic sequences} Given a model (thus, an equation)
\begin{equation}\label{eq:equation}
\mathcal E: y^2+a_1xy+a_3y=x^3+a_2x^2+a_4x
\end{equation}
for an elliptic curve containing the point $S=(0,0)$, and some point
$M$, denote the $x$-co-ordinate 
$x(M+hS)$ of $M+hS$ by $x_{M+hS}=-e_h$. Then straightforward
computations lead to several remarkable, and remarkably useful,
identities.

\begin{proposition}\label{pr:basic} There are constants $\alpha$,
$\beta$, $\gamma$, depending on the model $\mathcal E$ but
independent of the integer parameter
$h$, \emph{and of} `the translation' 
$M$, so that
\begin{gather}\label{eq:1}
e_{h-1}e_{h}^2e_{h+1}=\alpha^2 e_h-\beta\,;\\
\label{eq:2}
(e_{h-1}+e_{h+1})e_h^2=\gamma e_h-\alpha^2\,. 
\end{gather}
\end{proposition}
\noindent It is a straightforward exercise%
\footnote{It is worth noticing that \eqref{eq:2} follows from
\eqref{eq:1}. Indeed, we have 
$$
e_{h}e_{h+1}^2e_{h+2}-\alpha^2
e_{h+1}=e_{h-1}e_{h}^2e_{h+1}-\alpha^2 e_h\,.
$$
Dividing by $e_{h}e_{h+1}$ and cutely inserting
$e_{h}e_{h+1}$ on each side then yields
$$e_{h+1}e_{h+2}+\alpha^2/e_{h+1}+e_{h}e_{h+1}
=e_{h}e_{h+1}+\alpha^2/e_{h}+e_{h-1}e_{h}=\gamma\,,\quad\text{some
constant.}
$$}
to confirm such identities
by the formulaire for adding points on $\mathcal E$, 
see~\cite{SS}; the  arguments of
\cite{169}, making explicit the \cfe\ of Adams and Razar~\cite{AR},
provide a seemingly very different proof. We also mention
the following corollary.
\begin{corollary}\label{co:odd} Thus
$\alpha^2 (e_{h}+e_{h+1})=e_he_{h+1}(\gamma -e_he_{h+1})
+\beta$,
and therefore
\begin{equation}\label{eq:odd}
e_{h-1}e_{h}^2e_{h+1}^2e_{h+2}= \beta
e_he_{h+1}+(\alpha^4-\beta\gamma)\,.
\end{equation}
\end{corollary}
\begin{proof} Proposition~\ref{pr:basic} reports that 
$$
(X-e_{h-1}e_h)(X-e_he_{h+1})
=X^2-(\gamma e_h-\alpha^2)X/e_h+(\alpha^2e_h-\beta)\,;
$$
and then $X=e_he_{h+1}$ provides the `thus'. Therefore, indeed,
\begin{multline*}
e_{h-1}e_{h}^2e_{h+1}\cdot e_{h}e_{h+1}^2e_{h+2}=
(\alpha^2 e_h-\beta)(\alpha^2 e_{h+1}-\beta)\\
=\alpha^4 e_he_{h+1}-\alpha^2\beta(e_h+e_{h+1})+\beta^2
=e_he_{h+1}\bigl(\alpha^4+\beta(e_he_{h+1}-\gamma)\bigr)\,,
\end{multline*}
completing the proof.\end{proof} 

Further, define the \emph{elliptic sequence} $(A_h)$ by a pair of
initial values and the recursive definition 
\begin{equation}\label{eq:recursive}
A_{h-1}A_{h+1}=e_{h}A_h^2\,. 
\end{equation}
One checks readily that in immediate consequence
of \eqref{eq:recursive}:
\begin{multline}\label{eq:examples}
A_{h-2}A_{h+2}=e_{h-1}e_{h}^2e_{h+1}A_h^2\,;\quad\\ 
A_{h-1}A_{h+2}=e_{h}e_{h+1}A_hA_{h+1}\,;\quad 
A_{h-2}A_{h+3}=e_{h-1}e_{h}^2e_{h+1}^2e_{h+2}A_hA_{h+1}\,.
\end{multline}
In particular, multiplying \eqref{eq:1} by $A_h^2$, or \eqref{eq:odd}
by $A_hA_{h+1}$, yields the recursions
\begin{gather}\label{eq:recursion4}
A_{h-2}A_{h+2}=\alpha^2 A_{h-1}A_{h+1}-\beta A_h^2\,,\\
A_{h-2}A_{h+3}=\beta
A_{h-1}A_{h+2}+(\alpha^4-\beta\gamma)A_hA_{h+1}\,,
\label{eq:recursion5}\end{gather}
for the sequence $(A_h)$. So an {elliptic sequence} $(A_h)$ is
not quite the most general Somos~$4$ because the first coefficient
in the recursion is necessarily a square. 

\label{page:equivalence} However, replacing the
sequence
$(e_h)$ by $(\alpha e_h)$, thus the co-ordinates $x_{M+hS}$ by
$\alpha x_{M+hS}$,  transforms
$(A_h)$ into an
\emph{equivalent} sequence 
$(A'_h)$ with
$$
A'_h=\alpha^{h(h-1)/2}A_h\,,\quad\text{so that}\quad 
A'_{h-2}A'_{h+2}=\alpha^{-1} A'_{h-1}A'_{h+1}-\beta\alpha^{-4}
\smash{A'}_h^2\,.
$$
Thus any Somos~$4$ is, in the sense just described, at worst
equivalent to an elliptic sequence. Of course, in
place of our quadratic `twist' by $\alpha$, we could simply have
confessed to viewing the sequence as being an elliptic sequence over
a quadratic extension of the base field. 

Plainly, in the sequel we
may suppose without additional comment that claims we prove for
elliptic sequences hold appropriately for a general Somos~$4$
sequence.

\subsection{Singular elliptic sequences} Several of
our confident remarks above need to be announced more falteringly if 
some
$e_h$ vanishes. That case is $M+kS=\pm S$, some $k\in\Z$, and
thus,  by changing the translation
$M$ if necessary, there is no serious loss of generality in supposing
that in fact
$e_1=0$ and $W_0=0$. If 
some
$e_k$ vanishes, then $M+kS=\pm S$, some $k\in\Z$, and
thus,  by changing the translation
$M$ if necessary, there is no serious loss of generality in supposing
that in fact
$e_1=0$ and $W_0=0$. In this \emph{singular}\footnote{We make no
attempt to stick to Ward's terminolgy~\cite{Wa}; there, all the
sequences are `singular', in our sense of the adjective as
`special'.
 In \cite{Wa}  a sequence is called `singular' if it arises from
a singular cubic curve: such singularity is no issue for us.} case we set
$\e_h=-x_{hS}$ and define $W_{h-1}W_{h+1}=\e_hW_{h}^2$. So
$W_0=0$, and we may take
$W_1=1$; by \eqref{eq:1} we have $\e_2=\beta/\alpha^2$ and it is
reasonable to select $W_2=\alpha$, hence $W_3=\beta$; and --- we
leave the computation of $\e_3$ for the energetic reader (see
\cite{Sw} or~\cite{169}) ---
$W_4=-W_2^5+W_2W_3\gamma$. Of course, \eqref{eq:recursion4} becomes
\begin{equation}\label{eq:even'}
W_1^2A_{h-2}A_{h+2}=W_{2}^2A_{h-1}A_{h+1}-W_{1}W_{3}A_{h}^2\,;
\end{equation}
and \eqref{eq:recursion5} plainly is
\begin{equation}\label{eq:odd'}
W_1W_2A_{h-2}A_{h+3}=W_{2}W_3A_{h-1}A_{h+2}-W_{1}W_{4}A_{h}A_{h+1}\,.
\end{equation}
Here we have pretended to forget that $W_1=1$, the more vividly to
emphasise the anti-symmetry of the two-sided sequence $(W_h)$ and
related pattern. Of course $(W_h)$ is precisely the, well let's say
it, `untranslated' elliptic sequence discussed by Morgan
Ward~\cite{Wa}.

\subsection{Asides} Christine Swart~\cite{Sw} shows that the $A_h^2$ `try to be'
the denominators of the $x$~co-ordinates $-e_h=x_{M+hS}$  in that
they succeed in so being at worst up to finitely many primes
involved in the initial values and the defining recursion of the
sequence (and thus in the coefficients of the model $\mathcal E$ of
the underlying elliptic curve). More specifically, in the singular
case, Rachel Shipsey~\cite{Shi} confirms that if the model~$\mathcal
E$ is minimal integral with $\gcd(a_3, a_4)=1$ then $W_2\Div W_4$
guarantees 
$(W_h)$ is an exact  division sequence:
$\gcd (W_i,W_j)=W_{\gcd(i,j)}$.

If both $-x_S=\e_1=0$ \emph{and} $\e_{m+1}=0$, then the
sequence
$(\e_h)$ is periodic of period~$m$ --- for this case see~\cite{163}
and remarks at \cite[\S VIII]{Sw}--- but the singular elliptic
sequence
$(W_h)$ need be no more than quasi-periodic of quasi-period~$m$. We
skirt by the fact that then
$\e_0$, 
$\e_m$, $\ldots\,$ are infinite (so, of course $W_0$, $W_m$,  $\ldots\,$ must all
vanish) by noting in particular that the recursion relations for $(A_h)$ and
$(W_h)$ allow one to skip over and then fill in any difficulties;
we define the `undefined' portions of our sequences accordingly.

\subsection{Induction and symmetry} A surprisingly simple inductive
argument together with pleasing applications of symmetry suffice to
prove our main result: A Somos~$4$ also is a three-term Somos~$k$
for  $k=5$, $6$, $7$
$\ldots\,$.
\begin{theorem}\label{eq:main}
If, for $h\in\Z$, $A_{h-1}A_{h+1}=e_hA_{h}^2$ and
$W_{h-1}W_{h+1}=\e_hW_{h}^2$ then for $m\in\Z$
\begin{equation}\label{eq:evenm}
W_1^2A_{h-m}A_{h+m}=W_{m}^2A_{h-1}A_{h+1}-
W_{m-1}W_{m+1}A_{h}^2
\end{equation}
and
\begin{equation}\label{eq:oddm}
W_{1}W_{2}A_{h-m}A_{h+m+1}=W_{m}W_{m+1}A_{h-1}A_{h+2}-
W_{m-1}W_{m+2}A_{h}A_{h+1}\,.
\end{equation}
\end{theorem}
\begin{proof} We note that \eqref{eq:evenm} is
$$
\frac{A_{h-m}A_{h+m}}{A_h^2}=W_m^2
\left(\frac{A_{h-1}A_{h+1}}{A_h^2}
-\frac{W_{m-1}W_{m+1}}{W_m^2}\right)=W_m^2(e_h-\e_m)\,,
$$
so that, seeing that \eqref{eq:evenm} 
is trivially true for $m=0$ and
$m=1$, it suffices to show that
$$
\frac{A_{h-m-1}A_{h+m+1}}{A_h^2}=W_{m+1}^2(e_h-\e_{m+1})
$$
follows. However, by appropriate inductive hypotheses,
\begin{multline*}
\frac{A_{h-m-1}A_{h+m+1}}{A_h^2}\\=
\frac{A_{h-1-m}A_{h-1+m}}{A_{h-1}^2}\cdot
\frac{A_{h-1}A_{h+1}}{A_h^2}
\cdot\frac{A_{h+1-m}A_{h+1+m}}{A_{h+1}^2}
\Bigg/
\frac{A_{h-(m-1)}A_{h+(m-1)}}{A_h^2}
\\=
\frac{W_m^2}{W_{m-1}^2}\frac{(e_{h-1}-\e_{m})e_h^2(e_{h+1}-\e_{m})}
{(e_{h}-\e_{m-1})}.
\end{multline*}
This is in fact $W_{m+1}^2(e_h-\e_{m+1})$, as hoped for,  if and
only if
\begin{equation}\label{eq:symmetriceven}
(e_{h-1}-\e_{m})e_h^2(e_{h+1}-\e_{m})=
(\e_{m-1}-e_{h})\e_m^{\,2}(\e_{m+1}-e_{h})\,.
\end{equation}
We would like to be able to declare that \eqref{eq:symmetriceven} is
blatantly true by a principle of symmetry but, sadly,
a slightly more brutal argument seems necessary.

By Proposition~\ref{pr:basic} we have
\begin{multline*}
(e_{h-1}-\e_{m})e_h^2(e_{h+1}-\e_{m})\\
=e_{h-1}e_{h}^2e_{h+1}-\e_m(e_{h-1}+e_{h+1})e_h^2+\e_m^{\,2}e_h^2\\
=\alpha^2e_h-\beta-(\gamma\,\e_me_h-\alpha^2\,\e_m)+\e_m^{\,2}e_h^2\,,
\end{multline*}
making manifest the symmetry  $e_h\longleftrightarrow \e_m$.

Similarly, proving \eqref{eq:oddm} requires we show as inductive
step that
$$
\frac{W_1W_2}{W_{m+1}W_{m+2}}\frac{A_{h-(m+1)}A_{h+(m+2)}}{A_hA_{h+1}}
=(e_he_{h+1}-\e_{m+1}\e_{m+2}).
$$
But by permissible inductive hypotheses
\begin{multline*}
\frac{W_1W_2}{W_{m+1}W_{m+2}}\frac{A_{h-(m+1)}A_{h+(m+2)}}{A_hA_{h+1}}\\=
\frac{W_1W_2}{W_{m}W_{m+1}}\frac{A_{h-1-m}A_{h+m}}{A_{h-1}A_{h}}
\cdot\frac{W_m^2}{W_{m-1}W_{m+1}}\frac{W_{m+1}^2}{W_{m}W_{m+2}}
\frac{A_{h-1}A_{h+2}}{A_hA_{h+1}}\cdot\\\cdot
\frac{W_1W_2}{W_{m}W_{m+1}}\frac{A_{h-m+1}A_{h+m+2}}{A_hA_{h+1}}
\Bigg/
\frac{W_1W_2}{W_{m-1}W_{m}}\frac{A_{h-(m-1)}A_{h+m}}{A_hA_{h+1}}
\\=
\frac{(e_{h-1}e_h-\e_{m}\e_{m+1})e_he_{h+1}(e_{h+1}e_{h+2}-\e_{m}\e_{m+1})}
{\e_{m}\e_{m+1}(e_{h}e_{h+1}-\e_{m-1}\e_m)}.
\end{multline*}
So \eqref{eq:oddm} follows if and only if
\begin{multline*}
(e_{h-1}e_h-\e_{m}\e_{m+1})e_he_{h+1}(e_{h+1}e_{h+2}-\e_{m}\e_{m+1})\\
=(\e_{m-1}\e_m-e_{h}e_{h+1})\e_{m}\e_{m+1}(\e_{m+1}\e_{m+2}-e_{h}e_{h+1}).
\end{multline*}
Thus, here too, it suffices to notice that
\begin{multline*}
(e_{h-1}e_h-\e_{m}\e_{m+1})e_he_{h+1}(e_{h+1}e_{h+2}-\e_{m}\e_{m+1})\\
=e_{h-1}e_{h}^2e_{h+1}^2e_{h+2}-(e_{h-1}e_h+e_{h+1}e_{h+2})e_he_{h+1}\e_m\e_{m+1}
+e_he_{h+1}\e_m^{\,2}\e_{m+1}^{\,2}
\\=\beta
e_he_{h+1}+(\alpha^4-\beta\gamma)-\bigl(\alpha^2(e_h+e_{h+1})-2\beta\bigr)\e_m\e_{m+1} 
+e_he_{h+1}\e_m^{\,2}\e_{m+1}^{\,2}
\\=\beta e_he_{h+1}+(\alpha^4-\beta\gamma)-
\bigl(e_he_{h+1}(\gamma-e_he_{h+1})-\beta\bigr)\e_m\e_{m+1} 
+e_he_{h+1}\e_m^{\,2}\e_{m+1}^{\,2}
\end{multline*}
is visibly symmetric for $e_h\longleftrightarrow \e_m$.
\end{proof} 

\subsection{Comments}\subsubsection{}\label{ss:5}
Almost precisely the  argument just given shows also that: \emph{A
Somos~$5$ also is a three-term Somos~$k$ for  $k=7$, $9$, $11$
$\ldots\,$.}
Namely, sequences $(e_h)$ and $(c_h)$ give rise to a Somos~$5$
sequence $(B_h)$ by way of the definition
$c_hB_{h-1}B_{h+1}=e_hB_h^2$ for $h\in\Z$. Indeed, that yields
$$
c_{h-1}c_{h}^2c_{h+1}^2c_{h+2}W_1W_2B_{h-2}B_{h+3}
=c_{h}c_{h+1}W_{2}W_{3}B_{h-1}B_{h+2}
-W_{1}W_{4}B_{h}B_{h+1}\,;\
$$
and this relation has constant coefficients (thus, independent
of~$h$) exactly when $c_hc_{h+1}=:v$, say, is constant: so when
the sqeuences $(c_{2h})$ and$(c_{2h+1})$ are constant.

It now suffices to replace $e_h$ by $e_h/c_h$ in the argument above
(while not changing $\e_h$) to see that the just stated relation
implies that for all integers $h$ and $m$
\begin{multline}\label{eq:5}
v^{\frac12m(m+1)}W_1W_2B_{h-m}B_{h+m+1}
\\=vW_{m}W_{m+1}B_{h+1}B_{h+2}-W_{m-1}W_{m+2}B_{h}B_{h+1}\,.
\end{multline}
For example \cite{169}, $5$-Somos comes from the points $M+hS$ on the
elliptic curve 
$$
y^2+xy+6y=x^3+7x^2+12x\,,\qquad\text{with $M=(-2,-2)$ and
$S=(0,0)$\,.}
$$
In this case $W_1=1$, $W_2=6$, $W_3=6^2$, $W_4=-6^4$ so the choice
$v=6$ (more precisely, $c_0=2$, $c_1=3$) is felicitous.

\subsubsection{} Plainly
\begin{multline*}
A_{h-m}A_{h+m}W_n^2
=(W_{m}^2A_{h-1}A_{h+1}-W_{m-1}W_{m+1}A_{h}^2)W_n^2\\
= (W_{n}^2A_{h-1}A_{h+1}-W_{n-1}W_{n+1}A_{h}^2)W_m^2
-(W_{n}^2W_{m-1}W_{m+1}-W_{n-1}W_{n+1}W_{m}^2)A_h^2\\
=A_{h-n}A_{h+n}W_{m}^2-W_{m-n}W_{m+n}A_{h}^2\,,
\end{multline*}
confirming also that \eqref{eq:general} and \eqref{eq:redundant} are
indeed equivalent. Just so,
\begin{multline*}
A_{h-m}A_{h+m+1}W_nW_{n+1}
\\=(W_{m}W_{m+1}A_{h-1}A_{h+2}-W_{m-1}W_{m+2}A_{h}A_{h+1})W_nW_{n+1}/W_2\\
=(W_nW_{n+1}A_{h-1}A_{h+2}-W_{n-1}W_{n+2}A_{h}A_{h+1})W_{m}W_{m+1}/W_2\phantom{pushleft}
\\-(W_nW_{n+1}W_{m-1}W_{m+1}-W_{n-1}W_{n+2}W_{m}W_{m+1})A_hA_{h+1}/W_2\\
=W_{m}W_{m+1}A_{h-n}A_{h+n+1}-W_{m-n}W_{m+n+1}A_{h}A_{h+1}\,.
\end{multline*}

\subsubsection{} It warrants remark that elliptic curves play
at most an implicit role in our arguments. It suffices to start from
the identities of Proposition~\ref{pr:basic}. Moreover, \eqref{eq:1}
\emph{is} just the Somos~$4$ relation (if necessary by enlarging the
base field to include $\alpha$); and, as remarked, \eqref{eq:2}
follows.

\subsubsection{} Nonetheless, it is clear that all Somos~4 and
Somos~$5$ sequences \emph{are} elliptic sequences. We confine
ourselves here to just a sketch. Indeed, to be given a Somos~$4$
sequence is to be provided with a sequence
$(e_h)$ by way of
$e_h=A_{h-1}A_{h+1}/A_h^2$ where the $e_h$ satisfy \eqref{eq:1} and
hence all of Proposition~\ref{pr:basic}. Our main argument then
provides all the recurrence relations
$A_{h-m}A_{h+m}=W_{m}^2A_{h-1}A_{h+1}
-W_{m-1}W_{m+1}A_{h}^2$.
Those yield the sequence $(W_m)$
 --- which, see for example \cite{Shi} or \cite{Wa},
amounts  to having the related elliptic curve. Of course\footnote{Of
course this may not seem all that evident; fortunately, it is part
of the content of \cite{169} and~of~\cite{SS}.}
$e_0$ provides the `translation' $M$.

The Somos~$5$ case is a little less straightforward. Here  we have
$B_{h-1}B_{h+1}/B_h^2=e_h/c_h=:f_h$, say, where $c_hc_{h+1}=v$ is
constant --- equivalently, that
$(c_{2h})$ and
$(c_{2h+1})$ 	are constant sequences --- and we might define an
`equivalent' sequence $(A_h)$ by 
$$
B_{2h}=c_0^{-h(h+1)}c_1^{-h^2}A_{2h} \quad\text{and}\quad
B_{2h+1}=c_0^{-(h+1)^2}c_1^{-h(h+1)}A_{2h+1}\,. 
$$ 
Then $(A_h)$ is a Somos~$4$. Moreover, by arguments similar to those
we use above, if
\begin{multline*}\label{eq:}
A_{h-2}A_{h+2}=W_{2}^2A_{h-1}A_{h+1}-W_{3}A_{h}^2\,,
\\
\quad\text{then}\quad
A_{hk-2k}A_{hk+2k}=(W_{2k}/W_k)^2A_{hk-k}A_{hk+k}
     -(W_{3k}/W_k)A_{hk}^2
\end{multline*}
for all $k=1$, $2$, $3$, $\ldots\,$. Thus $(A_{hk})$ is a Somos~$4$
for all those $k$ and in particular it follows readily that both
$(B_{2h})$ and
$(B_{2h+1})$ are Somos~$4$ sequences.

\bibliographystyle{amsalpha}

\label{page:lastpage}
\end{document}